\documentclass[12pt,twoside]{article}
\usepackage{graphicx}
\usepackage{amsfonts}
\usepackage{amsbsy}
\usepackage{amssymb}
\usepackage{amsmath,tikz,pgfplots}
\usetikzlibrary{decorations.pathmorphing}
\usepackage{cite}
\usepackage{graphics}
\def\bpsp{\begin{pspicture}}
\def\epsp{\end{pspicture}}

\newcommand{\R}{\mathbb{R}}

\newtheorem{theorem}{Theorem}[section]
\newtheorem{remark}[theorem]{Remark}
\newtheorem{example}[theorem]{Example}
\newtheorem{lemma}[theorem]{Lemma}
\newtheorem{corollary}[theorem]{Corollary}
\newtheorem{definition}[theorem]{Definition}
\newtheorem{proposition}[theorem]{Proposition}
\newtheorem{note}{Note}
\newtheorem{case}{Case}
\newtheorem{conjecture}{Conjecture}
\newtheorem{question}{Question}

\newcommand{\bea}{\begin{eqnarray}}
\newcommand{\eea}{\end{eqnarray}}
\newcommand{\beq}{\begin{eqnarray*}}
\newcommand{\eeq}{\end{eqnarray*}}

\def\m4{\mbox{\rm ~(mod $4$)}}

\def \bd{\begin{definition}}
\def \ed{\end{definition}}
\def \bqu{\begin{question}}
\def \equ{\end{question}}
\def \bcc{\begin{conjecture}}
\def \ecc{\end{conjecture}}
\def \bt{\begin{theorem}}
\def \et{\end{theorem}}
\def \bl{\begin{lemma}}
\def \el{\end{lemma}}
\def \bc{\begin{corollary}}
\def \ec{\end{corollary}}
\def \be{\begin{equation}}
\def \ee{\end{equation}}
\def \ben{\begin{enumerate}}
\def \een{\end{enumerate}}
\def \ba{\begin{array}}
\def \ea{\end{array}}
\def \bp{\begin{proposition}}
\def \ep{\end{proposition}}
\def \bx{\begin{example}}
\def \ex{\end{example}}
\def \br{\begin{remark}}
\def \er{\end{remark}}
\def \bdsc{\begin{description}}
\def \edsc{\end{description}}

\def \bn{\begin{case}}
\def \en{\end{case}}
\def \bnt{\begin{note}}
\def \ent{\end{note}}
\def\1{1\!\!1}

\def\mm2{\mbox{\rm ~(mod $2$)}}
\def\m4{\mbox{\rm ~(mod $4$)}}

\def\qed{\nolinebreak\hfill\rule{.2cm}{.2cm}\par\addvspace{.5cm}}

\def\m{\mu}

\def\1{\textbf{1}}
\def\0{\textbf{0}}

\begin{document}
\title{On the eigenvalues of signed complete bipartite graphs }
\author{ S. Pirzada$^{a}$, Tahir Shamsher$^{b}$, Mushtaq A. Bhat$^{c}$ \\
$^{}$ {\em  $^{}$Department of Mathematics, University of Kashmir, Srinagar, Kashmir, India}\\
$^{}${\em  Department of Mathematics, National Institute of Technology, Srinagar, India}\\
$^a$pirzadasd@kashmiruniversity.ac.in,~~$^b$tahir.maths.uok@gmail.com,~~ \\$^c$mushtaqab@nitsri.net }
\date{}

\pagestyle{myheadings} \markboth{ Pirzada, Tahir, Bhat}{On  the eigenvalues of signed complete bipartite graphs }
\maketitle
\vskip 5mm
\noindent{\footnotesize \bf Abstract.} Let $\Gamma=(G,\sigma)$ be a signed graph, where $\sigma$ is the sign function on the edges of $G$. The adjacency matrix of $\Gamma=(G, \sigma)$ is a square matrix $A(\Gamma)=A(G, \sigma)=\left(a_{i j}^{\sigma}\right)$, where $a_{i j}^{\sigma}=\sigma\left(v_{i} v_{j}\right) a_{i j}$. In this paper, we determine the eigenvalues of the signed complete bipartite graphs.   Let $(K_{p, q},\sigma)$, $p\leq q$, be a signed complete bipartite graph with bipartition $(U_p, V_q)$, where $U_p=\{u_1,u_2,\ldots,u_p\}$ and $V_q=\{v_1,v_2,\ldots,v_q\}$. Let  $(K_{p, q},\sigma)[U_r\cup V_s]$, $r\leq p$ and $s\leq q $, be an induced signed subgraph on minimum vertices $r+s$,  which contains all negative edges of the signed graph $(K_{p, q},\sigma)$. We show that the multiplicity of  eigenvalue $0$ in $(K_{p, q},\sigma)$ is at least $ p+q-2k-2$, where $k=min(r,s)$. We determine the spectrum of signed complete bipartite graph whose negative edges induce disjoint complete bipartite subgraphs and path. We obtain  the spectrum of signed complete bipartite graph whose negative edges (positive edges) induce an $r-$ regular subgraph $H$. We find a relation between the eigenvalues of this  signed complete bipartite graph and the non-negative  eigenvalues of $H$.

\vskip 3mm



\section{Introduction}\label{sec1}
 A signed graph (or briefly sigraph) $\Gamma$ is an ordered pair $(G, \sigma)$, where $G=(V(G), E(G))$ is a graph (called the underlying graph), and $\sigma: E(G) \longrightarrow\{-1,1\}$ is a sign function defined on the edge set of $G$.  A signed graph is all-positive (all-negative) if all of its edges are positive (negative) and is denoted by $\Gamma=(G,+)$ (resp. $\Gamma=(G,-))$. The sign of a cycle in a signed graph is the product of the signs of its edges. A signed graph is said to be balanced if each of its cycle is positive, otherwise unbalanced.\\
\indent Let $A(G)=\left(a_{i j}\right)$ be the adjacency matrix of $G$. The adjacency matrix of a signed graph $\Gamma=(G, \sigma)$ is a square matrix $A(\Gamma)=A(G, \sigma)=\left(a_{i j}^{\sigma}\right)$, where $a_{i j}^{\sigma}=\sigma\left(v_{i} v_{j}\right) a_{i j}$. For a matrix $Z$, the characteristic polynomial $|x I-Z|$ will be denoted by $\phi(Z,x)$. If $\Gamma$ is a signed graph, we use $\phi(\Gamma,x)$ instead of $\phi(A(\Gamma),x)$. The eigenvalues of $A(\Gamma)$ are  the eigenvalues of the signed graph $\Gamma$. The set of all eigenvalues of $\Gamma$ along with their multiplicities is called the spectrum of the signed graph $\Gamma$. If the distinct eigenvalues of $\Gamma$ are $\mu_{1}>\dots>\mu_{k}$, and their multiplicities are $m\left(\mu_{1}\right), \ldots, m\left(\mu_{k}\right)$, then we write
$$
{\bf \sigma}(\Gamma)=\left(\begin{array}{ccc}
\mu_{1} & \ldots & \mu_{k} \\
m\left(\mu_{1}\right) & \ldots & m\left(\mu_{k}\right)
\end{array}\right).
$$
Two signed graphs $\Gamma_{1}=\left(G_{1}, \sigma_{1}\right)$ and $\Gamma_{2}=\left(G_{2}, \sigma_{2}\right)$ are isomorphic if there is a graph isomorphism $f: G_{1} \rightarrow G_{2}$ that preserves signs of the edges. If $\theta: V(G) \rightarrow\{+1,-1\}$ is the switching function, then switching of the signed graph $\Gamma=(G, \sigma)$ by $\theta$ means changing $\sigma$ to $\sigma^{\theta}$ defined by
$$
\sigma^{\theta}(u v)=\theta(u) \sigma(u v) \theta(v).
$$
For more information about switching, we refer to \cite{mb}.\\
Infact,  we observe that the sign function for the signed  subgraph   is the restriction of the signed graph $\Gamma$. For $X \subseteq V(G)$,  $\Gamma[X]$ denotes the  induced signed  subgraph formed by $X$, while $\Gamma-X=\Gamma[V(G) \backslash X].$ Sometimes, we also write $\Gamma-\Gamma[X]$ instead of $\Gamma-X$. Let $\left(G, K^{-}\right)$ ($\left(G, K^{+}\right)$) be the signed graph whose negative edges (positive edges) induce a subgraph $K$. As usual,  $K_{n}$ denotes the complete graph of order $n$. The complete bipartite graph  with two parts $U_p=\{u_1,u_2,\ldots,u_p\}$ and $V_q=\{v_1,v_2,\ldots,v_q\}$ as a partition of its vertex set is denoted by $K_{p, q}$. Also, $P_n$ denotes the path on $n$ vertices. $J_{r \times s}$ will denote an all-one matrix of size $r\times s$ and $O_{r\times s}$ will denote an all-zero matrix of size $r\times s$.\\

The rest  of the paper is organized as follows. In Section $2$, we give some preliminary results which will be used in the sequel. In Section $3$,  we show that the multiplicity of an eigenvalue $0$ in $(K_{p, q},\sigma)$ is at least $ p+q-2k-2$, where $k=min(r,s)$ and  $(K_{p, q},\sigma)[U_r\cup V_s]$, $r\leq p$ and $s\leq q $, is an induced signed subgraph on minimum vertices  $r+s$,  which contain all negative edges  of the signed graph $(K_{p, q},\sigma)$.  In section $4$, we determine the spectrum of the signed complete bipartite graph whose negative edges (positive edges) induce (i) disjoint complete bipartite subgraphs and (ii) a path. In Section $5$,
we determine the spectrum of the signed complete bipartite graph whose negative edges (positive edges) induce an $r-$ regular subgraph $H$. Also,  we obtain a relation between the eigenvalues of this  signed complete bipartite graph and the non-negative  eigenvalues of $H$.

\section{Preliminaries}\label{sec2}

 Let $l, m$ and $n$ be real numbers such that $l+m=n$. Then the two numbers $l$ and $m$ are symmetric with respect to $\frac{n}{2}$. Consider $\mu_{1},\mu_{2}, \ldots, \mu_{n}$ as the eigenvalues of the signed graph $\Gamma$. If for each $i$ there exists some $j$ such that $\mu_{i}+\mu_{j}=0$, then we say that the spectrum is symmetric with respect to $0$. It is well known that a graph which contains at least one edge is bipartite if and only if its spectrum considered as a set of points on the real axis is symmetric with respect to the origin. There exist nonbipartite signed graphs with this property as can be seen in \cite{mb1}. The following result can be seen in \cite{mb}.
 
 \begin{lemma} {\em\cite{mb}} \label{2.1}
 Let $\Gamma$ be a signed graph of order $n$. Then the following statements are equivalent.\\
(i) Spectrum of $\Gamma$ is symmetric about the origin,\\
(ii) $\phi(\Gamma,x)=x^{n}+\sum_{k=1}^{\left\lfloor\frac{n}{2}\right\rfloor}(-1)^{k} c_{2 k} x^{n-2 k}$, where $c_{2 k}$ are non negative integers for all $k=1,2, \ldots,\left\lfloor\frac{n}{2}\right\rfloor$,\\
(iii) $\Gamma$ and $-\Gamma$ are cospectral, where $-\Gamma$ is the signed graph obtained by negating sign of each edge of $\Gamma$.
 \end{lemma}

Consider a matrix $M$ having the block form as follows.
\begin{equation}
M=\left(\begin{array}{ccccc}
A & \beta & \cdots & \beta & \beta \\
\beta^{\top} & B & \cdots & C & C \\
\vdots & \vdots & \cdots & \vdots & \vdots \\
\beta^{\top} & C & \cdots & B & C \\
\beta^{\top} & C & \cdots & C & B
\end{array}\right)
\end{equation}
where $A \in R^{t \times t}, ~\beta \in R^{t \times s}$ and $B, C \in R^{s \times s}$, such that $n=t+c s$, with $c$ being the number of copies of $B$. The spectrum of this matrix can be obtained as the union of the spectrum of smaller matrices using the following technique given in \cite{fs}. In the statement of the following theorem, $\sigma^{(k)}(Z)$ denotes the multi-set formed by $k$ copies of the spectrum of $Z$, denoted by $\sigma(Z)$.

\begin{lemma}\label{2.2} Let $M$ be a matrix of the form given in $(2.1)$ with $c \geq 1$ copies of the block $B$. Then\\
(i) $\sigma(B-C) \subseteq \sigma(M)$ with multiplicity $c-1$,\\
(ii) $\sigma(M) \backslash \sigma^{(c-1)}(B-C)=\sigma\left(M^{\prime}\right)$ is the set of the remaining $t+s$ eigenvalues of $M$, where
$$
M^{\prime}=\left(\begin{array}{cc}
A & \sqrt{c} \cdot \beta \\
\sqrt{c} \cdot \beta^{\top} & B+(c-1) C
\end{array}\right).
$$
\end{lemma}

\begin{lemma}\label{2.3}
Let $
X=\left(\begin{array}{cc}
O_{p\times p} & A_{p\times q} \\
 A_{q\times p}^{\top} & O_{q\times q}
\end{array}\right)
$ be a real symmetric square matrix of order $p + q$, $q\geq p$. Then \\
(i) $m(0)\geq q-p$, \\
(ii) $\pm \sqrt{\mu} \in \sigma(X)$, where $\mu$ is an eigenvalue of a positive semidefinite square matrix $A_{p\times q}A_{q\times p}^{\top} $.
\end{lemma}
{\bf Proof.} By Schur complement formula, the determinant  of a $2\times 2$ block matrix  is given by
$$
\left|\begin{array}{ll}
A & B \\
C & D
\end{array}\right|=\left|D \| A-B D^{-1} C\right|,
$$
where $A$ and $D$ are square blocks and $D$ is nonsingular. So, we have
$$
\begin{aligned}
\phi(X,x) &=\left|\begin{array}{cc}
x I_{p} & -A_{p\times q} \\
-A_{q\times p}^{\top} & x I_{q}
\end{array}\right|=x^{q}\left|\left(x I_{p}\right)-A_{p\times q}\left(x I_{q}\right)^{-1} A_{q\times p}^{\top}\right|=x^{q}\left|\frac{1}{x}\left(x^{2} I_{p}-A_{p\times q}A_{q\times p}^{\top}\right)\right| \\
&= x^{q-p} \phi\left(A_{p\times q}A_{q\times p}^{\top}, x^{2}\right).
\end{aligned}
$$
This completes the proof. \qed

\begin{corollary} \label{2.4}
Let $
X=\left(\begin{array}{cc}
O_{p\times p} & A_{p\times p} \\
 A_{p\times p } & O_{p\times p}
\end{array}\right)
$ be a real symmetric square matrix of order $2p$. Then $\pm \mu \in \sigma(X)$, where $\mu$ is an eigenvalue of the square matrix $A_{p\times p} $.
\end{corollary}

{\bf Remark 2.1} Let $(K_{p, q},\sigma)$ be a signed complete bipartite graph  with bipartition $(U_p, V_q)$, where $U_p=\{u_1,u_2,\ldots,u_p\}$ and $V_q=\{v_1,v_2,\ldots,v_q\}$. Then with a suitable labelling of the vertices of $(K_{p, q},\sigma)$,  its adjacency matrix is given by  $$ A(K_{p, q},\sigma)
=\left(\begin{array}{cc}
O_{p\times p} & B_{p\times q} \\
 B_{q\times p } ^{\top}& O_{q\times q}
\end{array}\right).
$$
In view of  Lemma \ref{2.3}, we observe that the spectrum of $(K_{p, q},\sigma)$ is related with the spectrum of the matrix $B_{p\times q}B_{q\times p}^{\top}$. Thus from here onwards, we focus on the matrix $B_{p\times q}$ and we call it as the spectral block of the adjacency matrix of the signed graph $(K_{p, q},\sigma)$.
 
\section{Multiplicity of the eigenvalue $0$ in $(K_{p, q},\sigma)$}\label{sec3}

 In this section, we obtain a lower bound for the multiplicity of the eigenvalue $0$ in $\Gamma=(K_{p, q},\sigma)$ for any sign function $\sigma$, subject to the condition.
 
 \begin{theorem}\label{3.1}
 Let $(K_{p, q},\sigma)$, $p\leq q$, be a signed complete bipartite graph and let $(K_{p, q},\sigma)[U_r\cup V_s]$, $r\leq p$ and $s\leq q $, be its induced signed subgraph on minimum vertices $r+s$,  which contains all negative edges of the signed graph $(K_{p, q},\sigma)$. Then $m(0)\geq p+q-2k-2$, where $k=min(r,s)$.
 \end{theorem}
 {\bf Proof.} Note that the order of  $(K_{p, q},\sigma)[U_r\cup V_s]$ is $r+s$.   With a suitable labelling of the vertices of $(K_{p, q},\sigma)$,  the adjacency matrix is given by  $$ A(K_{p, q},\sigma)
=\left(\begin{array}{cc}
O_{p\times p} & B_{p\times q} \\
 B_{q\times p } ^{\top}& O_{q\times q}
\end{array}\right),
$$
 where, $B_{p \times q}$ is the spectral block of the adjacency matrix of the signed graph $(K_{p, q},\sigma)$. By Lemma \ref{2.3}, we get
 \begin{equation}
\begin{aligned}
\phi(A(K_{p, q},\sigma),x) = x^{q-p} \phi\left(B_{p\times q}B_{q\times p}^{\top}, x^{2}\right).
\end{aligned}
\end{equation}

As $(K_{p, q},\sigma)[U_r\cup V_s]$ is an induced signed subgraph on minimum vertices $r+s$,  which contain all negative edges of the signed graph $(K_{p, q},\sigma)$, we have
$$ B_{p \times q}
=\left(\begin{array}{cc}
X_{r\times s} & J_{r\times q-s} \\
 J_{p-r\times s } & J_{p-r\times q-s}
\end{array}\right),
$$
where $X_{r \times s}$ is the spectral block of the adjacency matrix of the signed graph $(K_{p, q},\sigma)[U_r\cup V_s]$. The transpose of a $2\times 2$ block matrix  is given by $$\left(\begin{array}{ll}A & B \\ C & D\end{array}\right)^{\top}=\left(\begin{array}{ll}A^{\top} & C^{\top} \\ B^{\top} & D^{\top}\end{array}\right).$$
Together  with the fact that $J_{m\times n}J_{n\times m} = n J_{m \times m}$, this yields
$$
\begin{aligned}
 B_{p \times q}B_{q \times p}^{\top}
&=\left(\begin{array}{cc}
X_{r\times s} & J_{r\times q-s} \\
 J_{p-r\times s } & J_{p-r\times q-s}
\end{array}\right)\times \left(\begin{array}{cc}
X_{s\times r}^{\top} & J_{s\times p-r} \\
 J_{q-s\times r } & J_{q-s\times p-r}
\end{array}\right)\\
&=\left(\begin{array}{cc}
X_{r\times s}X_{s\times r}^{\top}+(q-s)J_{r \times r} & X_{r\times s}J_{s\times p-r}+(q-s)J_{r\times p-r} \\
  J_{p-r\times s}X_{s\times r}^{\top}+(q-s)J_{p-r\times r} & sJ_{p-r\times p-r}+(q-s)J_{p-r\times p-r}
\end{array}\right)\\
&=\left(\begin{array}{cc}
X_{r\times s}X_{s\times r}^{\top}+(q-s)J_{r \times r} & X_{r\times s}J_{s\times p-r}+(q-s)J_{r\times p-r} \\
  J_{p-r\times s}X_{s\times r}^{\top}+(q-s)J_{p-r\times r} & qJ_{p-r\times p-r}
\end{array}\right).
\end{aligned}
$$
Now, it is easy to see that $X_{r\times s}J_{s\times 1}+(q-s)J_{r\times 1}=Y+(q-s)J_{r\times 1} $, where $Y$ is the column vector of the row sums of the matrix $X_{r\times s}$. Let $Z=[Y+(q-s)J_{r\times 1} ~~Y+(q-s)J_{r\times 1} ~~\cdots~~ Y+(q-s)J_{r\times 1}]\in R^{r\times p-r}$ be a matrix of order $r\times p-r$. Then, we have
\begin{equation}
\begin{aligned}
 B_{p \times q}B_{q \times p}^{\top}
=\left(\begin{array}{cc}
X_{r\times s}X_{s\times r}^{\top}+(q-s)J_{r \times r} & Z \\
  Z^{\top} & qJ_{p-r\times p-r}
\end{array}\right).
\end{aligned}
\end{equation}
The matrix $B_{p \times q}B_{q \times p}^{\top}$ has a special kind of symmetry. Taking $A=X_{r\times s}X_{s\times r}^{\top}+(q-s)J_{r \times r}$, $\beta = Y+(q-s)J_{r\times 1}$, $B=[q]$ and $C=[q]$ in $(2.1)$,  from Lemma \ref{2.2}, we get $\sigma^{p-r-1}(B-C)=\sigma^{p-r-1}([0])\subseteq \sigma(B_{p \times q}B_{q \times p}^{\top})$. Again by Eq. $(3.2)$, Eq. $(3.3)$ and Lemma \ref{2.2}, we obtain
\begin{equation}
\phi(A(K_{p, q},\sigma),x) = x^{\alpha} \phi\left(Z_1, x^{2}\right),
\end{equation}
where $\alpha=q+p-2r-2$ and $Z_1=\left(\begin{array}{cc}
X_{r\times s}X_{s\times r}^{\top}+(q-s)J_{r \times r} & \sqrt{p-r}(Y+(q-s)J_{r\times 1}) \\
  \sqrt{p-r}(Y+(q-s)J_{r\times 1})^{\top} & q(p-r)
\end{array}\right)$.\\
Also, we have
$$
\begin{aligned}
 B_{q \times p}^{\top}B_{p \times q}
=\left(\begin{array}{cc}
X_{s\times r}^{\top}X_{r\times s}+(p-r)J_{s \times s} & X_{s\times r}^{\top}J_{r\times q-s}+(p-r)J_{s\times q-s} \\
  J_{q-s\times r}X_{r\times s}+(p-r)J_{q-s\times s} & pJ_{q-s\times q-s}
\end{array}\right).
\end{aligned}
$$
Now,  $X_{s\times r}^{\top}J_{r\times 1}+(p-r)J_{s\times 1}=Y^{\prime}+(p-r)J_{s\times 1} $, where $Y^{\prime}$ is the column vector of the column sums of the matrix $X_{r\times s}$. Let $Z^{\prime}=[Y^{\prime}+(p-r)J_{s\times 1} ~~Y^{\prime}+(p-r)J_{s\times 1} ~~\cdots~~ Y^{\prime}+(p-r)J_{s\times 1}]\in R^{s\times q-s}$ be a matrix of order $s\times q-s$. Then,
\begin{equation}
\begin{aligned}
 B_{q \times p}^{\top}B_{p \times q}
=\left(\begin{array}{cc}
X_{s\times r}^{\top}X_{r\times s}+(p-r)J_{s \times s} & Z^{\prime} \\
  Z'^{\top} & pJ_{q-s\times q-s}
\end{array}\right).
\end{aligned}
\end{equation}
 Taking $A=X_{s\times r}^{\top}X_{r\times s}+(p-r)J_{s \times s}$, $\beta = Y'+(p-r)J_{s\times 1}$, $B=[p]$ and $C=[p]$ in $(2.1)$,  from Lemma \ref{2.2}, we get $\sigma^{q-s-1}(B-C)=\sigma^{q-s-1}([0])\subseteq {\bf \sigma}(B_{p \times q}^{\top}B_{p \times q})$. Note that the eigenvalues of $ B_{q \times p}^{\top}B_{p \times q}$ are given by the eigenvalues of $ B_{p \times q}B_{q \times p}^{\top}$, together with the eigenvalue $0$ of multiplicity $q-p$. Therefore, by Eq. $(3.2)$, Eq. $(3.5)$ and Lemma \ref{2.2}, we obtain
\begin{equation}
\phi(A(K_{p, q},\sigma),x) = x^{\zeta} \phi\left(Z_2, x^{2}\right),
\end{equation}
where $\zeta=q+p-2s-2$ and $Z_2=\left(\begin{array}{cc}
X_{s\times r}^{\top}X_{r\times s}+(p-r)J_{s \times s} & \sqrt{q-s}(Y'+(p-r)J_{s\times 1}) \\
  \sqrt{q-s}(Y'+(p-r)J_{s\times 1})^{\top} & p(q-s)
\end{array}\right)$. Hence the result follows by Eq. $(3.4)$ and Eq. $(3.6)$.\qed
As $(K_{p, q},\sigma)$ is a signed bipartite  graph and therefore its spectrum is symmetric about the origin. Thus, the following is an immediate consequence of Theorem \ref{3.1} and Lemma \ref{2.1}.
\begin{figure}
\centering
	\includegraphics[scale=.8]{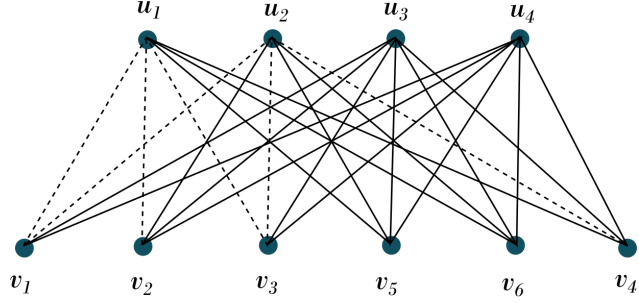}
	\caption{The signed graph $(K_{4,6},\sigma)$.}
	\label{Figure 1}
\end{figure}
\begin{corollary}\label{3.2}
 Let $(K_{p, q},\sigma)$, $p\leq q$, be a signed complete bipartite graph and let $(K_{p, q},\sigma)(U_r\cup V_s)$, $r\leq p$ and $s\leq q $, be its induced subgraph on minimum vertices $r+s$,  which contains all  positive edges of the signed graph $(K_{p, q},\sigma)$. Then $m(0)\geq p+q-2k-2$, where $k=min(r,s)$.
\end{corollary}

\noindent {\bf Example 3.1} Consider a signed complete bipartite graph $(K_{4,6},\sigma)$ as shown in Figure 1. Plain lines denote the positive edges and dashed lines denote the negative edges. It contains an induced signed subgraph $(K_{4,6},\sigma)[U_2,V_4]$ on $6$ vertices which contain all negative edges of $(K_{4,6},\sigma)$.  Here, we have $p=4$, $q=6$ and $k=2$. Therefore, by Theorem $3.1$,  $m(0)\geq 4$. The spectral block of the adjacency matrix of the induced signed subgraph  $(K_{4,6},\sigma)[U_2,V_4]$ is given as
 $$X_{2\times4}=\left(\begin{array}{cccc}
-1 & -1&-1 &1\\
-1&1&-1&-1
\end{array}\right).$$
Therefore, by Eq. $(3.4)$, we get
 $$\phi(A(K_{4, 6},\sigma),x) = x^{4} \phi\left(\left(\begin{array}{ccc}
6 & 2&0 \\6&2&0\\
  0&0 & 12
\end{array}\right), x^{2}\right).$$ Thus, it is easy to see that
$$
{\bf \sigma}((K_{4, 6},\sigma))=\left(\begin{array}{ccccccc}
2\sqrt{3} & 2\sqrt{2}& 2 & 0 & -2 & -2\sqrt{2} & -2\sqrt{3} \\
1 & 1 & 1 & 4 & 1& 1& 1
\end{array}\right).
$$

\section{Eigenvalues of $\left(K_{p,q}, K_{r, s}^{-}\right)$}\label{sec4}

 We begin this section with a signed complete bipartite graph $\left(K_{p,q}, K_{r, s}^{-}\right)$ whose negative edges induce a subgraph $K_{r, s}$.
 
\begin{theorem}\label{4.1}
 Let $\left(K_{p,q}, K_{r, s}^{-}\right)$, $p \leq q$,  $r\leq p$ and $s\leq q$, be a signed complete bipartite graph whose negative edges induce a subgraph  $ K_{r,s}$ of order $r+s$. Then the spectrum of $\left(K_{p,q}, K_{r, s}^{-}\right)$ is given as
 $$
{\bf \sigma}(\left(K_{p,q}, K_{r, s}^{-}\right))=\left(\begin{array}{ccccc}
\mu_1 & \mu_2 &  0 & -\mu_2 & -\mu_1 \\
1 & 1 & p+q-4 & 1& 1\\
\end{array}\right),
$$
where $$ \mu_1, \mu_2 = \sqrt{\frac{rq+q(p-r) \pm \sqrt{r(q^2r+2(p-r)(q-2s)^2)+q^2(p-r)^2}}{2}}.$$
\end{theorem}
{\bf Proof.} By Eq. $(3.4)$, we have
\begin{equation} \phi\left((K_{p,q}, K_{r, s}^{-}\right),x) = x^{\alpha} \phi\left(\left(\begin{array}{cc}
X_{r\times s}X_{s\times r}^{\top}+(q-s)J_{r \times r} & \sqrt{p-r}(Y+(q-s)J_{r\times 1}) \\
  \sqrt{p-r}(Y+(q-s)J_{r\times 1})^{\top} & q(p-r)
\end{array}\right), x^{2}\right),
\end{equation}
where $\alpha=q+p-2r-2$ and $Y$ is the column vector of the row sums of spectral block $X_{r\times s}$ of the adjacency matrix of an induced signed subgraph $K_{r, s}$, whose all edges are negative.  Clearly, $X_{r\times s}X_{s\times r}^{\top}+(q-s)J_{r \times r}=-J_{r \times s}\times -J_{s \times r}+(q-s)J_{r \times r}=qJ_{r \times r} $ and  $Y+(q-s)J_{r\times 1}=(q-2s)J_{r\times 1}$. Therefore,  Eq.$(4.7)$ takes the form
\begin{equation} \phi\left((K_{p,q}, K_{r, s}^{-}\right),x) = x^{\alpha} \phi\left(\left(\begin{array}{cc}
qJ_{r \times r} & \sqrt{p-r}(q-2s)J_{r\times 1} \\
  \sqrt{p-r}(q-2s)J_{1\times r} & q(p-r)
\end{array}\right), x^{2}\right).
\end{equation}
It can be easily  seen that the real symmetric matrix $$Z_1=\left(\begin{array}{cc}
qJ_{r \times r} & \sqrt{p-r}(q-2s)J_{r\times 1} \\
  \sqrt{p-r}(q-2s)J_{1\times r} & q(p-r)
\end{array}\right)$$ has rank $2$. Now, let $x_1$ and $x_2$ be the non zero eigenvalues of $Z_1$. We have\\
\begin{equation}
x_1+x_2=tr(Z_1) =rq+q(p-r).
\end{equation}
Also,
\begin{equation}
x_1^2+x_2^2=tr(Z_1^2) =r(q^2r+2(p-r)(q-2s)^2)+q^2(p-r)^2.
\end{equation}
Eqs. $(4.9)$ and $(4.10)$, imply that
\begin{equation}
x_1,~x_2 =\frac{rq+q(p-r) \pm \sqrt{r(q^2r+2(p-r)(q-2s)^2)+q^2(p-r)^2}}{2}.
\end{equation}
Thus, Eq. $(4.8)$ yields that
\begin{equation*} \phi\left((K_{p,q}, K_{r, s}^{-}\right),x) = x^{p+q-4}(x^4-(x_1+x_2)x^2+x_1x_2),
\end{equation*}
where $x_1$ and $x_2$ are given in Equation $(4.11)$. This proves the result. \qed

\begin{corollary}\label{4.2}
 Let $\left(K_{p,q}, K_{r, s}^{+}\right)$, $p \leq q$, $r\leq p$ and $s\leq q$, be a signed complete bipartite graph whose positive edges induce a subgraph  $ K_{r,s}$ of order $r+s$. Then the spectrum of $\left(K_{p,q}, K_{r, s}^{+}\right)$ is given as
 $$
{\bf \sigma}(\left(K_{p,q}, K_{r, s}^{+}\right))=\left(\begin{array}{ccccc}
\mu_1 & \mu_2 &  0 & -\mu_2 & -\mu_1 \\
1 & 1 & p+q-4 & 1& 1\\
\end{array}\right),
$$
where $$ \mu_1, \mu_2 = \sqrt{\frac{rq+q(p-r) \pm \sqrt{r(q^2r+2(p-r)(q-2s)^2)+q^2(p-r)^2}}{2}}.$$
\end{corollary}

Now, we consider the signed complete bipartite graph $\left(K_{p,q}, \sigma\right)$ whose negative edges form  disjoint  subgraphs   $K_{r, s}$ of different orders.
\begin{theorem}\label{4.3}
 Let $\left(K_{p,q},\sigma\right)$, $p \leq q$, be a signed complete bipartite graph whose negative edges induce disjoint complete bipartite subgraphs  of different orders $r_1+s_1$, $r_2+s_2$, $\cdots$, $r_k+s_k$ such that $\sum\limits_{i=1}^{k}r_i=r$,  $\sum\limits_{i=1}^{k}s_i=s$, $r\leq p$ and $s\leq q$. Then the characteristic polynomial of $\left(K_{p,q}, \sigma\right)$ is given as
 $$
\phi\left((K_{p,q}, \sigma\right),x) = x^{p+q-2k-2} \phi\left(Z', x^{2}\right),$$
where $$Z'=
\left(\begin{array}{ccccc}
r_1c_{11} & r_2c_{12}& \cdots & r_kc_{1k} &c(q-2s_1)\\
  r_1c_{21} & r_2c_{22} &\cdots & r_kc_{2k} & c(q-2s_2)
  \\
  \vdots & \vdots & \ddots & \vdots &\vdots\\
  r_1c_{k1} & r_2c_{k2} & \cdots & r_kc_{kk} &c(q-2s_k)\\
  r_1c(q-2s_1) & r_2c(q-2s_2) &\cdots & r_kc(q-2s_k) & q(p-r)
\end{array}\right)$$ is a positive semidefnite matrix of order $k+1$, $c=\sqrt{p-r}$, $c_{ij}=q$ if $i=j$ and $c_{ij}=q-2s_i-2s_j$ otherwise.
\end{theorem}
{\bf Proof.} Consider the matrix given in Eq. $(3.4)$
\begin{equation}
Z_1=
\left(\begin{array}{cc}
X_{r\times s}X_{s\times r}^{\top}+(q-s)J_{r \times r} & \sqrt{p-r}(Y+(q-s)J_{r\times 1}) \\
  \sqrt{p-r}(Y+(q-s)J_{r\times 1})^{\top} & q(p-r)
\end{array}\right),
\end{equation}
where $Y$ is the column vector of the row sums of the spectral block $X_{r\times s}$ of the adjacency matrix of an induced signed subgraph of $\left(K_{p,q},\sigma\right)$ which contains all its  negative edges. Hence, with a suitable relabelling of vertices of the induced signed subgraph, we have
$$X_{r\times s}=
\left(\begin{array}{cccc}
-J_{r_1 \times s_1} & J_{r_1 \times s_2}& \cdots &J_{r_1 \times s_k}\\
  J_{r_2 \times s_1} & -J_{r_2 \times s_2} &\cdots & J_{r_2 \times s_k}\\
  \vdots & \vdots & \ddots & \vdots\\
  J_{r_k \times s_1} & J_{r_k \times s_2} & \cdots & -J_{r_k \times s_k}
\end{array}\right),$$
where $J_{r_i \times s_i}$ is the spectral block of the adjacency matrix of the complete bipartite subgraph $K_{r_i,s_i}$, $i=1, 2, \ldots, k$, $\sum\limits_{i=1}^{k}r_i=r$ and $\sum\limits_{i=1}^{k}s_i=s$. Now, it is easy to obtain
$$X_{r\times s}X_{s\times r}^{\top}=
\left(\begin{array}{cccc}
b_{11}J_{r_1 \times r_1} & b_{12}J_{r_1 \times r_2}& \cdots & b_{1k}J_{r_1 \times r_k}\\
  b_{21}J_{r_2 \times r_1} & b_{22}J_{r_2 \times r_2} &\cdots & b_{2k}J_{r_2 \times r_k}\\
  \vdots & \vdots & \ddots & \vdots\\
  b_{k1}J_{r_k \times r_1} & b_{k2}J_{r_k \times r_2} & \cdots & b_{kk}J_{r_k \times r_k}
\end{array}\right),$$
where, $b_{ij}=\sum\limits_{i=1}^{k}s_i=s$ if $i=j$ and $b_{ij}=s-2s_i-2s_j$ otherwise. As $Y$ is the column vector of the row sums of the spectral block $X_{r\times s}$, therefore the matrix $Z_1$ given in $(4.12)$ takes the form
$$Z_1=
\left(\begin{array}{ccccc}
c_{11}J_{r_1 \times r_1} & c_{12}J_{r_1 \times r_2}& \cdots & c_{1k}J_{r_1 \times r_k} &c(q-2s_1)J_{r_1 \times 1}\\
  c_{21}J_{r_2 \times r_1} & c_{22}J_{r_2 \times r_2} &\cdots & c_{2k}J_{r_2 \times r_k} & c(q-2s_2)J_{r_2 \times 1}\\
  \vdots & \vdots & \ddots & \vdots &\vdots\\
  c_{k1}J_{r_k \times r_1} & c_{k2}J_{r_k \times r_2} & \cdots & c_{kk}J_{r_k \times r_k} &c(q-2s_k)J_{r_k \times 1}\\
  c(q-2s_1)J_{1 \times r_1} & c(q-2s_2)J_{1 \times r_2} &\cdots & c(q-2s_k)J_{1 \times r_k} & q(p-r)J_{1 \times 1}
\end{array}\right),$$
where $c=\sqrt{p-r}$, $c_{ij}=q$ if $i=j$ and $c_{ij}=q-2s_i-2s_j$ otherwise. Clearly, the matrix $Z_1$ has equitable quotient matrix $Z'$, where
$$Z'=
\left(\begin{array}{ccccc}
r_1c_{11} & r_2c_{12}& \cdots & r_kc_{1k} &c(q-2s_1)\\
  r_1c_{21} & r_2c_{22} &\cdots & r_kc_{2k} & c(q-2s_2)
  \\
  \vdots & \vdots & \ddots & \vdots &\vdots\\
  r_1c_{k1} & r_2c_{k2} & \cdots & r_kc_{kk} &c(q-2s_k)\\
  r_1c(q-2s_1) & r_2c(q-2s_2) &\cdots & r_kc(q-2s_k) & q(p-r)
\end{array}\right).$$ Now by [Theorem $3.1$,\cite{q}], ${\bf \sigma}(Z_1)={\bf \sigma}(Z')\cup \left(\begin{array}{c}
0\\r-k
\end{array}\right)$, where $Z'$ is equitable quotient matrix of $Z_1$ and is given as
$$Z'=
\left(\begin{array}{ccccc}
r_1c_{11} & r_2c_{12}& \cdots & r_kc_{1k} &c(q-2s_1)\\
  r_1c_{21} & r_2c_{22} &\cdots & r_kc_{2k} & c(q-2s_2)
  \\
  \vdots & \vdots & \ddots & \vdots &\vdots\\
  r_1c_{k1} & r_2c_{k2} & \cdots & r_kc_{kk} &c(q-2s_k)\\
  r_1c(q-2s_1) & r_2c(q-2s_2) &\cdots & r_kc(q-2s_k) & q(p-r)
\end{array}\right),$$
\begin{figure}
\centering
	\includegraphics[scale=.8]{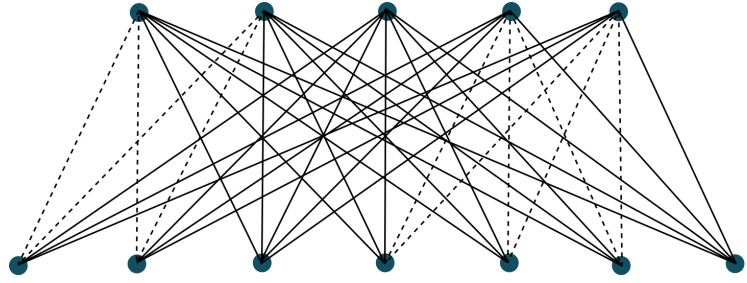}
	\caption{Signed graph  whose negative edges induce two disjoint complete bipartite subgraphs.}
	\label{Figure 1}
\end{figure}
where $c=\sqrt{p-r}$, $c_{ij}=q$ if $i=j$ and $c_{ij}=q-2s_i-2s_j$ otherwise. As ${\bf \sigma}(Z_1)={\bf \sigma}(Z')\cup \left(\begin{array}{c}
0\\r-k
\end{array}\right)$, therefore the result follows by Eq. $(3.4)$ and Eq. $(4.12)$. \qed

\noindent{\bf Example 4.1} Consider the signed complete bipartite graph $(K_{5,7},\sigma)$ as shown in Figure 2.  Here, we have $p=5$, $q=7$, $r_1=2$, $s_1=2$, $r_2=2$, $s_2=3$, $r=r_1+r_2=4$ and $s=s_1+s_2=5$. Therefore, by Theorem \ref{4.3},  we get
 $$\phi(A(K_{5, 7},\sigma),x) = x^{6} \phi\left(\left(\begin{array}{ccc}
14 & -6&3 \\-6&14&1\\
  6&2 & 7
\end{array}\right), x^{2}\right).$$ Thus, it is easy to see that
$$
{\bf \sigma}((K_{5, 7},\sigma))=\left(\begin{array}{ccccccc}
4.50 & 3.37& 1.82 & 0 & -1.82 & -3.37 & -4.50 \\
1 & 1 & 1 & 6 & 1& 1& 1
\end{array}\right).
$$
\begin{corollary}\label{4.4}
 Let $\left(K_{p,q},\sigma\right)$, $p \leq q$, be a signed complete bipartite graph whose positive edges induce disjoint complete bipartite subgraphs  of different orders $r_1+s_1$, $r_2+s_2$, $\cdots$, $r_k+s_k$ such that $\sum\limits_{i=1}^{k}r_i=r$, $\sum\limits_{i=1}^{k}s_i=s$, $r\leq p$ and $s\leq q$. Then the characteristic polynomial of $\left(K_{p,q}, \sigma\right)$ is given as
 $$
\phi\left((K_{p,q}, \sigma\right),x) = x^{p+q-2k-2} \phi\left(Z', x^{2}\right),$$
where $$Z'=
\left(\begin{array}{ccccc}
r_1c_{11} & r_2c_{12}& \cdots & r_kc_{1k} &c(q-2s_1)\\
  r_1c_{21} & r_2c_{22} &\cdots & r_kc_{2k} & c(q-2s_2)
  \\
  \vdots & \vdots & \ddots & \vdots &\vdots\\
  r_1c_{k1} & r_2c_{k2} & \cdots & r_kc_{kk} &c(q-2s_k)\\
  r_1c(q-2s_1) & r_2c(q-2s_2) &\cdots & r_kc(q-2s_k) & q(p-r)
\end{array}\right)$$ is a positive semidefnite matrix of order $k+1$, $c=\sqrt{p-r}$, $c_{ij}=q$ if $i=j$ and $c_{ij}=q-2s_i-2s_j$ otherwise.
\end{corollary}
We conclude this section with the following result whose proof can be obtained in a similar way as in Theorem \ref{4.3}.
\begin{theorem}\label{4.5}
 $(i)$ Let $\left(K_{p,q},P_{2r}^-\right)$, $p \leq q$ and $r\geq 1$, be a signed complete bipartite graph whose negative edges induce a path on $2r$ vertices. Then the characteristic polynomial of $\left(K_{p,q}, P_{2r}^-\right)$ is given as
 $$
\phi\left((K_{p,q}, P_{2r}^-\right),x) = x^{p+q-2r-2} \phi\left(Z', x^{2}\right),$$
where
\begin{equation*}
Z'=
\left(\begin{array}{cccccccc}
q & q-2 & q-6 & q-6 & \cdots & q-6&c(q-2)\\
 q-2 &q &q-4 & q-8 &\cdots & q-8&c(q-4)\\
q-6 &q-4 &q &q-4&\ddots &\vdots &\vdots\\
q-6 &q-8 &q-4 &q &\ddots & q-8&c(q-4)\\
\vdots &\vdots &\ddots &\ddots  &\ddots &q-4&c(q-4)\\
q-6 &q-8 &\hdots &q-8 &q-4 &q&c(q-4)\\
c(q-2)&c(q-4)&\hdots&c(q-4)&c(q-4)&c(q-4)&q(p-r)

\end{array}\right)
\end{equation*}
 is a positive semidefnite matrix of order $r+1$ and $c=\sqrt{p-r}$.\\
  $(ii)$ Let $\left(K_{p,q},P_{2r+1}^-\right)$, $p \leq q$ and $r\geq 1$, be a signed complete bipartite graph whose negative edges induce a path on $2r+1$ vertices with both pendent vertices of the path $P_{2r+1}$ in $U_p$. Then the characteristic polynomial of $\left(K_{p,q}, P_{2r+1}^-\right)$ is given as
 $$
\phi\left((K_{p,q}, P_{2r+1}^-\right),x) = x^{p+q-2r-2} \phi\left(Z', x^{2}\right),$$
where
\begin{equation*}
Z'=
\left(\begin{array}{ccccccccc}
q & q-2 & q-6 & q-6 & \cdots & q-6&q-4&c(q-2)\\
 q-2 &q &q-4 & q-8 &\cdots & q-8&q-6&c(q-4)\\
q-6 &q-4 &q &q-4&\ddots &\vdots &\vdots&\vdots\\
q-6 &q-8 &q-4 &q &\ddots & q-8&q-6 &c(q-4)\\
\vdots &\vdots &\ddots &\ddots  &\ddots &q-4&q-6 &c(q-4)\\
q-6 &q-8 &\hdots &q-8 &q-4 &q& q-2&c(q-4)\\
q-4&q-6&\hdots&q-6& q-6&q-2 &q &c(q-2)\\
c(q-2)&c(q-4)&\hdots&c(q-4)&c(q-4)&c(q-4)&c(q-2) &q(p-r)

\end{array}\right)
\end{equation*}
 is a positive semidefnite matrix of order $r+2$ and $c=\sqrt{p-r}$.\\
 $(iii)$ Let $\left(K_{p,q},P_{2r+1}^-\right)$, $p \leq q$ and $r\geq 1$, be a signed complete bipartite graph whose negative edges induce a path on $2r+1$ vertices with both pendent vertices of the path $P_{2r+1}$ in $V_q$. Then the characteristic polynomial of $\left(K_{p,q}, P_{2r+1}^-\right)$ is given as
 $$
\phi\left((K_{p,q}, P_{2r+1}^-\right),x) = x^{p+q-2r-4} \phi\left(Z', x^{2}\right),$$
where
\begin{equation*}
Z'=
\left(\begin{array}{cccccccc}
q & q-4 & q-8 & q-8 & \cdots & q-8&c(q-4)\\
 q-4 &q &q-4 & q-8 &\cdots & q-8&c(q-4)\\
q-8 &q-4 &q &q-4&\ddots &\vdots &\vdots\\
q-8 &q-8 &q-4 &q &\ddots & q-8&c(q-4)\\
\vdots &\vdots &\ddots &\ddots  &\ddots &q-4&c(q-4)\\
q-8 &q-8 &\hdots &q-8 &q-4 &q&c(q-4)\\
c(q-4)&c(q-4)&\hdots&c(q-4)&c(q-4)&c(q-4)&q(p-r)

\end{array}\right)
\end{equation*}
 is a positive semidefnite matrix of order $r+1$ and $c=\sqrt{p-r}$.
\end{theorem}
{\bf Example 4.2} Let $\left(K_{p,q},P_{5}^-\right)$, $p \leq q$, be a signed complete bipartite graph whose negative edges induce a path on $5$ vertices with both pendent vertices of the path $P_{5}$ in $V_q$.  By Theorem \ref{4.5} (part $(iii)$), the characteristic polynomial of $\left(K_{p,q}, P_{5}^-\right)$ is given by
 $$
\phi\left((K_{p,q}, P_{5}^-\right),x) = x^{p+q-6} \phi\left(\left(\begin{array}{ccc}
q & q-4 &c(q-4)\\
 q-4 &q &c(q-4)\\
c(q-4)&c(q-4)&q(p-r)

\end{array}\right), x^{2}\right),$$
where $c=\sqrt{p-r}$. To determine the spectrum of $(K_{p,q}, P_{5}^-)$, it is enough to consider the matrix
$$Z'=\left(\begin{array}{ccc}
q & q-4 &c(q-4)\\
 q-4 &q &c(q-4)\\
c(q-4)&c(q-4)&q(p-r)

\end{array}\right).$$
Clearly, $4$ is an eigenvalue of the matriz $Z'$ corresponding to an eigenvector $(1,-1,0)^{\top}$. To compute the other two eigenvalues of $Z'$, we use the fact that the sum and product of the eigenvalues of $Z'$ are equal to the trace and determinant respectively. Then, we obtain the eigenvalues as
\begin{equation*}
 \frac{pq-4 \pm \sqrt{p^2q^2-56pq+128p+96q-240}}{2}.
\end{equation*}
Thus, the spectrum of $(K_{p,q}, P_{5}^-)$ is given as
$$
{\bf \sigma}((K_{p,q}, P_{5}^-))=\left(\begin{array}{ccccccc}
\mu_1 & \mu_2 &2&  0 &-2& -\mu_2 & -\mu_1 \\
1 & 1&1 & p+q-6 &1& 1& 1\\
\end{array}\right),
$$
where $$ \mu_1, \mu_2 = \sqrt{\frac{pq-4 \pm \sqrt{p^2q^2-56pq+128p+96q-240}}{2}}.$$

\section{Eigenvalues of $\left(K_{p,q}, H_{{\bf r}, n}^{-}\right)$}\label{sec4}

The signed complete bipartite graph $\Gamma$ whose negative edges induce a $1$-regular graph of different orders has been studied in \cite{ak}. In this section, we consider  signed complete bipartite graph $\left(K_{p,q}, H_{{\bf r}, n}^{-}\right)$, $p \leq q$, whose negative edges induce an ${\bf r}$-regular subgraph $H$ (not necessarily connected) of order $n$.  We find a relation between the eigenvalues of this signed complete bipartite graph and the non-negative eigenvalues of $H$. The other eigenvalues of $\left(K_{p,q}, H_{{\bf r}, n}^{-}\right)$ are also determined. We start with the following lemma.
\begin{lemma}\label{5.1}
Let $\left(K_{k,k}, H_{{\bf r}, 2k}^{-}\right)$ be a  signed complete bipartite graph whose negative edges induce an ${\bf r}$-regular subgraph $H$ of order $2k.$ If the eigenvalues of $H$ are $\mu_{1}={\bf r} \geq \mu_{2} \geq \cdots \geq \mu_{2k}=-{\bf r}$, then $-2 \mu_{i}$ is an eigenvalue of $\left(K_{k,k}, H_{{\bf r}, 2k}^{-}\right)$ for $i=2, \ldots, 2k-1.$ Moreover, the other two eigenvalues of $\left(K_{k,k}, H_{{\bf r}, 2k}^{-}\right)$ are $k-2{\bf r}$ and $-k+2{\bf r}$.
\end{lemma}
{\bf Proof.} Let $A(H,-)=-A(H)$ be the adjacency matrix of $(H,-).$  Therefore, with a suitable labelling of the vertices of $\left(K_{k,k}, H_{{\bf r}, 2k}^{-}\right)$, we observe that
\begin{equation}
A\left(K_{k,k}, H_{r, 2k}^{-}\right)=\left(\begin{array}{cc}
O_{k\times k} & A_{k\times k} \\
 A_{k\times k } & O_{k\times k}
\end{array}\right)=A(K_{k,k})-2A(H),
\end{equation}
where the  $(k-2{\bf r})$-regular symmetric matrix $A_{k\times k}$ is the spectral block of the adjacency matrix of the signed graph $\left(K_{k,k}, H_{{\bf r}, 2k}^{-}\right)$. As the matrices $A(K_{k,k})$ and $A(H)$  commute, therefore they are simultaneously diagonalizable. Let $\{x_{1}, x_{2}, \ldots, x_{2k}\}$ be an orthogonal basis of $\mathbb{R}^{2k}$ consisting of the eigenvectors of $A(H)$ and $A(K_{k,k})$ with $x_{1}=J_{2k\times1}=(1, \ldots, 1)^{T} \in \mathbb{R}^{2k}$. Then, we have
$$(A(K_{k,k})-2A(H))x_1=(k-2{\bf r})x_1.$$
Thus, $(k-2{\bf r})$ is an eigenvalue of $A(K_{k,k})-2A(H)$. To find the other eigenvalues of $A(K_{k,k})-2A(H)$, we use the facts that ${\bf \sigma}( A(K_{k,k})-2A(H))\subseteq {\bf \sigma}(A(K_{k,k}))+\sigma(-2A(H))$ and the spectrum of $\left(K_{k,k}, H_{{\bf r}, 2k}^{-}\right)$ is  symmetric with respect to origin. Thus,
$$(A(K_{k,k})-2A(H))x_i=-2\mu_ix_i, ~~i=2,3,\ldots,2k-1$$
and
$$(A(K_{k,k})-2A(H))x_{2k}=(-k+2{\bf r})x_{2k}.$$
This proves the result.\qed

\begin{theorem}
Let $\left(K_{p,q}, H_{{\bf r}, 2k}^{-}\right)$, $p \leq q$, be a  signed complete bipartite graph whose negative edges induce an ${\bf r}$-regular subgraph $H$ of order $2k.$ Then the following statements hold:\\
$(i)$ $m(0)\geq p+q-2k-2.$\\
$(ii)$ If the first $k$ largest non-negative eigenvalues of $H$ are $\mu_{1}={\bf r} \geq \mu_{2} \geq \dots \geq \mu_{k}\geq 0$, then $\pm ~2 \mu_{i}$ is an eigenvalue of $\left(K_{p,q}, H_{{\bf r}, 2k}^{-}\right)$ for $i=2, \ldots, k.$ Moreover, the other four eigenvalues of $\left(K_{p,q}, H_{{\bf r}, 2k}^{-}\right)$ are $$
\pm \sqrt{\frac{pq +(k-2{\bf r})^2-k^2 \pm \sqrt{(pq +(k-2{\bf r})^2-k^2)^2-4((k-2{\bf r})^2+k(q-k)-\frac{k(q-2{\bf r})^2}{q})(q(p-k))}}{2}}.
$$
\end{theorem}
{\bf Proof.}  Consider the matrix which is given in Eq. $(3.4)$
\begin{equation}
Z_1=
\left(\begin{array}{cc}
X_{r\times s}X_{s\times r}^{\top}+(q-s)J_{r \times r} & \sqrt{p-r}(Y+(q-s)J_{r\times 1}) \\
  \sqrt{p-r}(Y+(q-s)J_{r\times 1})^{\top} & q(p-r)
\end{array}\right),
\end{equation}
where, $Y$ is the column vector of the row sums of the spectral block $X_{r\times s}$ of the adjacency matrix of the induced signed subgraph  $\left(K_{k,k}, H_{{\bf r}, 2k}^{-}\right)$ of  $\left(K_{p,q}, H_{{\bf r}, 2k}^{-}\right)$ which contains all the  negative edges. By Eq. $(5.13)$, it is easy to see that $r=s=k$, $X_{r\times s}X_{s\times r}^{\top}=A_{k\times k}^2$ and $Y=(k-2{\bf r}) J_{k\times 1}$. Now, the matrix $Z_1$ takes the form
\begin{equation*}
Z_1=
\left(\begin{array}{cc}
A_{k\times k}^2+(q-k)J_{k \times k} & \sqrt{p-k}(q-2{\bf r})J_{k\times 1} \\
  \sqrt{p-k}(q-2{\bf r})J_{1\times k}^{\top} & q(p-k)
\end{array}\right).
\end{equation*}
 The matrix $A_{k\times k}^2$ is  $(k-2{\bf r})^2$-regular and hence commutes with $(q-k)J_{k \times k}$. Thus, it is easy to see that $(k-2{\bf r})^2-(q-k)$ is an eigenvalue of $A_{k\times k}^2+(q-k)J_{k \times k}$ corresponding to an eigenvector $J_{k\times1}$.  Also, by Eq. $(5.13)$, Corollary \ref{2.4} and Lemma \ref{5.1}, we have  $$(A_{k\times k}^2+(q-k)J_{k \times k})x_i=4\mu_{i}^2 x_i, i=2,\ldots, k,$$
 where $\{x_1,x_2,\ldots, x_k\}$ is an orthogonal basis of $\R^k$ with $x_1=J_{k\times1}$ and $\mu_i$ is non-negative eigenvalue of $H$. Define $y_i=[x_i~~ 0]^{\top}\in \R^{k+1}$, $i=2,\ldots, k$. Then $$Z_1y_i=4\mu_i^2y_i, ~i=2,\ldots,k.$$
 Therefore, $4\mu_i^2$, $i=2,\ldots, k$ is an eigenvalue of $Z_1$. Let $\alpha_1$ and $\alpha_2$  be the other two  eigenvalues of $Z_1$. We have\\
\begin{equation*}
\alpha_1+\alpha_2+\sum\limits_{i=2}^k4\mu_i^2=tr(Z_1) =k(q-2{\bf r})+q(p-k)
\end{equation*}
and
\begin{equation*}
(k-2{\bf r})^2+\sum\limits_{i=2}^k 4\mu_i^2= tr(A_{k\times k}^2) =k(k-2{\bf r}).
\end{equation*}
This yields that
\begin{equation}
\alpha_1+ \alpha_2 =pq +(k-2{\bf r})^2-k^2.
\end{equation}
By Schur complement formula, the determinant  of a $2\times 2$ block matrix $Z_1$ is given by
$$
\left|\begin{array}{cc}
A_{k\times k}^2+(q-k)J_{k \times k} & \sqrt{p-k}(q-2{\bf r})J_{k\times 1} \\
  \sqrt{p-k}(q-2{\bf r})J_{1\times k}^{\top} & q(p-k)
\end{array}\right|=\left|q(p-k)\right| \left|A_{k\times k}^2+(q-k)J_{k \times k}-\frac{(q-2{\bf r})^2}{q}J_{k\times k}\right|.
$$
Now, clearly the eigenvalues of the matrix $A_{k\times k}^2+(q-k)J_{k \times k}-\frac{(q-2{\bf r})^2}{q}J_{k\times k}$ are $(k-2r)^2+k(q-k)-\frac{k(q-2{\bf r})^2}{q}$ and $4\mu_i^2$, $i=2,\ldots, k$, where $\mu_i$ is the non-negative eigenvalue of $H$. Thus, we have
\begin{equation}
\alpha_1 \alpha_2= ((k-2{\bf r})^2+k(q-k)-\frac{k(q-2{\bf r})^2}{q})(q(p-k)).
\end{equation}
Equations $(5.15)$ and  $(5.16)$ imply that
\begin{equation*}
\alpha_1, \alpha_2= \frac{pq +(k-2{\bf r})^2-k^2 \pm \sqrt{(pq +(k-2{\bf r})^2-k^2)^2-4((k-2{\bf r})^2+k(q-k)-\frac{k(q-2{\bf r})^2}{q})(q(p-k))}}{2}.
\end{equation*}
Hence, by Eq. $(3.4)$, we have $m(0)\geq p+q-2k-2$ and with the fact that $\pm \alpha $ is an eigenvalue of  $\left(K_{p,q}, H_{{\bf r}, 2k}^{-}\right)$ whenever $\alpha^{2}$ is an eigenvalue of $Z_1$, the proof follows. \qed

\begin{corollary}
Let $\left(K_{p,q}, H_{{\bf r}, 2k}^{+}\right)$ be a  complete bipartite signed  graph whose positive edges induce an ${\bf r}$-regular subgraph $H$ of order $2k.$ Then the following statements hold:\\
$(i)$ $m(0)\geq p+q-2k-2.$\\
$(ii)$ If the first $k$ largest non-negative eigenvalues of $H$ are $\mu_{1}={\bf r} \geq \mu_{2} \geq \dots \geq \mu_{k}\geq 0$. Then $\pm ~2 \mu_{i}$ is an eigenvalue of $\left(K_{p,q}, H_{{\bf r}, 2k}^{+}\right)$ for $i=2, \ldots, k.$ Moreover, the other four  eigenvalues of $\left(K_{p,q}, H_{{\bf r}, 2k}^{+}\right)$ are $$
\pm \sqrt{\frac{pq +(k-2{\bf r})^2-k^2 \pm \sqrt{(pq +(k-2{\bf r})^2-k^2)^2-4((k-2{\bf r})^2+k(q-k)-\frac{k(q-2{\bf r})^2}{q})(q(p-k))}}{2}}.
$$
\end{corollary}

\begin{corollary}
Let $\left(K_{p,q}, H_{{\bf r}, 2k}^{}\right)$ be a  signed complete bipartite graph whose negative edges (positive edges) induce an ${\bf r}$-regular subgraph $H$ of order $2k.$ Then the signed graph $\left(K_{p,q}, H_{{\bf r}, 2k}^{}\right)$ is nonsingular if and only if the graph $H$ is nonsingular and $p=q=k$.
\end{corollary}

\noindent{\bf Acknowledgements.}   This research is supported by SERB-DST research project number CRG/2020/000109.  The research of Tahir Shamsher is supported by SRF financial assistance by the Council of Scientific and Industrial Research (CSIR), New Delhi, India.

\end{document}